\documentclass[11pt,twoside]{amsart}
\usepackage{amsxtra}
\usepackage{color}
\usepackage{amsopn}
\usepackage{amsmath,amsthm,amssymb,arydshln}
\usepackage{mathrsfs,mathtools}
\usepackage{stmaryrd}
\usepackage{hyperref}
\usepackage{enumerate}
\usepackage{xcolor}
\usepackage{pifont}
\usepackage{adjustbox}

\newtheorem{theorem}{Theorem}[section]
\newtheorem{corollary}[theorem]{Corollary}
\newtheorem{proposition}[theorem]{Proposition}

\newtheorem*{theorem*}{Theorem}
\theoremstyle{definition}
\newtheorem{definition}[theorem]{Definition}

\newtheorem{example}[theorem]{Example}
\newtheorem{remark}[theorem]{Remark}

\newcommand{\R}{{\mathbb R}}

\newcommand{\tu}{\tau_1}
\newcommand{\td}{\tau_2}

\newcommand{\W}{\wedge}

\newcommand{\f}{\varphi}
\newcommand{\fz}{\varphi_{\sst0}}

\newcommand{\SU}{{\rm SU}}

\newcommand{\G}{{\rm G}}
\newcommand{\K}{{\rm K}}

\newcommand{\ddt}{\frac{\partial}{\partial t}}

\newcommand{\Ric}{{\rm Ric}}

\newcommand{\Scal}{{\rm Scal}}

\newcommand{\fra}{\mathfrak{a}}

\newcommand{\frg}{\mathfrak{g}}

\newcommand{\frh}{\mathfrak{h}}

\newcommand{\frn}{\mathfrak{n}}

\newcommand{\frs}{\mathfrak{s}}

\newcommand{\st}{\ |\ }

\newcommand{\diag}{{\rm diag}}
\newcommand{\sst}{\scriptscriptstyle}

\numberwithin{equation}{section}

\title[On G$_2$-structures, special metrics and related flows]{On G$_{\mathbf2}$-structures, special metrics and related flows}
 
\author{Marisa Fern\'andez}
\address{Universidad del Pa\'{\i}s Vasco  (UPV / EHU), Facultad de Ciencia y Tecnolog\'{\i}a, 
Departamento de Matem\'aticas, Apartado 644, 48080 Bilbao, Spain}
\email{marisa.fernandez@ehu.es}
\author{Anna Fino} 
\address{Dipartimento di Matematica ``G.~Peano'' \\ Universit\`a degli Studi di Torino\\
Via Carlo Alberto 10\\
10123 Torino\\ Italy}
\email{annamaria.fino@unito.it}
\author{Alberto Raffero}
\address{Dipartimento di Matematica ``G.~Peano'' \\ Universit\`a degli Studi di Torino\\
Via Carlo Alberto 10\\
10123 Torino\\ Italy}
\email{alberto.raffero@unito.it}

\begin{document}
\begin{abstract} 
We review results about $\G_2$-structures in relation to the existence of special metrics, such as Einstein metrics and Ricci solitons,  
and the evolution under the Laplacian flow on non-compact homogeneous spaces. 
We also discuss some examples in detail.
\end{abstract}
\maketitle

%%%%%%%%%%%%%%%%%%%%%%%%%%%%%%%%%%%%%%%%%%%%%%%%%%%%%%%%%%%%%%%%%%%%%%%%%%%%%%%%%%%%%%%%%
%%%%%%%%%%%%%%%%%%%%%%%%%%%%%%%%%%%%%%%%%%%%%%%%%%%%%%%%%%%%%%%%%%%%%%%%%%%%%%%%%%%%%%%%%
%																INTRODUCTION 
%%%%%%%%%%%%%%%%%%%%%%%%%%%%%%%%%%%%%%%%%%%%%%%%%%%%%%%%%%%%%%%%%%%%%%%%%%%%%%%%%%%%%%%%%
%%%%%%%%%%%%%%%%%%%%%%%%%%%%%%%%%%%%%%%%%%%%%%%%%%%%%%%%%%%%%%%%%%%%%%%%%%%%%%%%%%%%%%%%%
\section{Introduction} 

A $\G_2$-structure on a seven-dimensional manifold $M$ is characterized by the existence of a globally defined 3-form $\f$ which can be pointwise written as 
\[
\f = e^{127} + e^{347}+ e^{567} + e^{135} - e^{146} - e^{236} - e^{245}, 
\]
with respect to a suitable basis $\{e^1, \ldots , e^7\}$ of the cotangent space. 
Here, the shorthand $e^{ijk}$ stands for $e^i\wedge e^j\wedge e^k$.
Such a 3-form $\f$ gives rise to a Riemannian metric $g_\f$ and an orientation $dV_\f$ on $M.$ 

The intrinsic torsion of a $\G_2$-structure $\f$ can be identified with the covariant derivative $\nabla^\f\f$, $\nabla^\f$ being the Levi Civita connection of $g_\f$. 
By \cite{FeGr}, it vanishes identically if and only if both $d\f=0$ and $d*_\f\f=0$, where $*_\f$ denotes the Hodge operator of $g_\f$. 
When this happens, the $\G_2$-structure is said to be {\em torsion-free}, its associated Riemannian metric $g_\f$ is Ricci-flat 
and the corresponding Riemannian holonomy group is a subgroup of the exceptional Lie group $\G_2$.

$\G_2$-structures can be divided into classes, which are characterized by the expression of the exterior derivatives $d\f$ and $d*_\f\f$ \cite{Bry,FeGr}.  
A $\G_2$-structure $\f$ is called {\em closed} (or {\em calibrated} according to \cite{HaLa}) if $d\f=0$, while it is called {\em coclosed} (or {\em cocalibrated}) if $d*_\f\f=0$.

Since the Ricci tensor and the scalar curvature of the metric induced by a $\G_2$-structure can be expressed in terms of the intrinsic torsion \cite{Bry}, 
it may happen that certain restrictions on the curvature give rise to some constraints on the intrinsic torsion. 
For instance, a calibrated $\G_2$-structure on a compact manifold induces an Einstein metric if and only if it is also cocalibrated, i.e., if and only if it is torsion-free \cite{Bry,ClIv}. 
A natural problem consists then in investigating whether this happens also in the non-compact case, and whether similar results also hold when the metric is a Ricci soliton.  
These problems were studied for calibrated $\G_2$-structures on homogeneous spaces in \cite{FFM,FFM2}, 
and for the wider class of locally conformal calibrated $\G_2$-structures in \cite{FiRa1}. 
We shall review the results in Section \ref{G2SpecMetrSect}. 

A useful tool to study geometric structures on manifolds is represented by geometric flows. 
Let $M$ be a 7-manifold endowed with a calibrated $\G_2$-structure $\fz$. 
The {\em Laplacian flow} starting from $\fz$ is the initial value problem
\[
\begin{cases}
\ddt \f(t) = \Delta_{\f(t)}\f(t),\\
d \f(t)=0,\\
\f(0)=\fz,
\end{cases}
\]
where $\Delta_\f\f=dd^*\f+d^*d\f$ is the Hodge Laplacian of $\f$ with respect to the metric $g_\f$. 
This flow was introduced by Bryant in \cite{Bry} to study 7-manifolds admitting calibrated $\G_2$-structures. 
Short-time existence and uniqueness of the solution when $M$ is compact were proved in the unpublished paper \cite{BrXu}. 
Recently, the analytic and geometric properties of the Laplacian flow have been deeply investigated in the series of papers \cite{LoWe1,LoWe2,LoWe3}. 
In particular, the authors obtained a long-time existence result, and they proved that the solution exists for all positive times and it 
converges to a torsion-free $\G_2$-structure modulo diffeomorphism 
provided that the initial datum $\fz$ is sufficiently close to a given torsion-free $\G_2$-structure. 

The first noncompact examples with long-time existence of the solution were obtained on seven-dimensional nilpotent Lie groups in \cite{FFM2},   
while further solutions on solvable Lie groups were described in \cite{FiRa,Lau2,Lau3,Nic}. 
Moreover, a cohomogeneity one solution converging to a torsion-free $\G_2$-structure on the 7-torus was worked out in \cite{HWY}.  
 In Section \ref{LapFlowSect}, we shall discuss the results on nilpotent Lie groups obtained in \cite{FFM2}.

%%%%%%%%%%%%%%%%%%%%%%%%%%%%%%%%%%%%%%%%%%%%%%%%%%%%%%%%%%%%%%%%%%%%%%%%%%%%%%%%%%%%%%%%%
%%%%%%%%%%%%%%%%%%%%%%%%%%%%%%%%%%%%%%%%%%%%%%%%%%%%%%%%%%%%%%%%%%%%%%%%%%%%%%%%%%%%%%%%%
%																PRELIMINARIES 
%%%%%%%%%%%%%%%%%%%%%%%%%%%%%%%%%%%%%%%%%%%%%%%%%%%%%%%%%%%%%%%%%%%%%%%%%%%%%%%%%%%%%%%%%
%%%%%%%%%%%%%%%%%%%%%%%%%%%%%%%%%%%%%%%%%%%%%%%%%%%%%%%%%%%%%%%%%%%%%%%%%%%%%%%%%%%%%%%%%
\section{Preliminaries}
Let $M$ be a seven-dimensional manifold endowed with a $\G_2$-structure $\f$. 
The Riemannian metric $g_\f$ and the volume form $dV_\f$ are determined by $\f$ via the equation
\[
g_\f(X,Y)\,dV_\f = \frac16\,\iota_X\f\W\iota_Y\f\W\f,
\]
for all vector fields $X,Y$ on $M.$

The vanishing of the intrinsic torsion $T_\f$  of a $\G_2$-structure $\f$ can be stated in the following equivalent ways.

\begin{theorem}[\cite{FeGr}] 
Let $\f$ be a $\G_2$-structure on a seven-dimensional manifold $M$. Then, the following conditions are equivalent:
\begin{enumerate}[a)]
\item the intrinsic torsion of $\f$ vanishes identically;
\item $\nabla^\f  \f=0$, where $\nabla^\f$ denotes the Levi Civita connection of $g_\f$;
\item $d\varphi=0$ and $d*_\f\f=0$;
\item $\mathrm{Hol}(g_{\f})$ is isomorphic to a subgroup of $\G_{2}$.
\end{enumerate}
\end{theorem}

A $\G_2$-structure satisfying any of the above conditions is said to be {\em torsion-free} or {\em parallel}.
By \cite{Bon}, the Riemannian metric induced by a torsion-free $\G_2$-structure $\f$ is Ricci-flat,  i.e., $\Ric(g_\f)=0$.

More generally, as the intrinsic torsion  $T_\f$ is a section of a vector bundle over $M$ with fibre ${\R^7}^*\otimes \frg_2^\perp$,   
$\G_2$-structures can be divided into classes according to the vanishing of the components of $T_\f$ with respect to the $\G_2$-irreducible decomposition
\[
{\R^7}^* \otimes {\mathfrak g}_2^\perp 	\cong 	{\mathcal X}_1\oplus {\mathcal X}_2\oplus {\mathcal X}_3\oplus {\mathcal X}_4 
								=		 \R \oplus {\mathfrak g}_2 \oplus S^2_0 (\R^7) \oplus \R^7,
\]
where $S^2_0 (\R^7)$ denotes the space of traceless symmetric 2-tensors and $\frg_2=\mathrm{Lie}(\G_2)$. 
This gives rise to sixteen classes of $\G_2$-structures, which were first described in \cite{FeGr}. 

By \cite{Bry}, it is also possible to characterize each class in terms of the exterior derivatives $d\f$ and $d*_\f\f$. 
In detail, the spaces $\Lambda^k({\R^7}^*)$, $k=2,3$, admit the following $\G_2$-irreducible decompositions (cf.~\cite{Bry87})
\[
\renewcommand\arraystretch{1.2}
\begin{array}{rcl}
\Lambda^2({\R^7}^*) &=& \Lambda^2_{7} ({\R^7}^*) \oplus \Lambda^2_{14} ({\R^7}^*),\\
\Lambda^3({\R^7}^*) &=& \Lambda^3_{1} ({\R^7}^*)\oplus \Lambda^3_7 ({\R^7}^*)\oplus \Lambda^3_{27}({\R^7}^*),
\end{array}
\renewcommand\arraystretch{1}
\]
where the subscript in $\Lambda^k_r ({\R^7}^*)$ denotes the dimension of the summand as an irreducible $\G_2$-module, 
and $\Lambda^2_{14} ({\R^7}^*)\cong\frg_2$, $\Lambda^3_{27} ({\R^7}^*)\cong S^2_0 (\R^7)$. 
Consequently, on $M$ there exist a unique function $\tau_0\in\mathcal{C}^\infty(M)$ and unique differential forms $\tau_1\in\Omega^1(M)$, 
$\tau_2\in\Omega^2_{14}(M)\coloneqq\left\{\alpha\in\Omega^2(M) \st \alpha\W*_{\f}\f = 0\right\}$, 
$\tau_3\in\Omega^3_{27}(M)\coloneqq \left\{\beta\in\Omega^3(M) \st \beta\W\f=0,~\beta\W*_\f\f=0 \right\}$ such that 
\begin{equation}\label{IntTorFor}
\begin{array}{rcl}
d \varphi &=& \tau_0 *_\f \varphi + 3 \tau_1 \wedge \varphi + *_\f \tau_3,
\\[4pt] d*_\f\varphi &=& 4 \tau_1 \wedge *_\f \varphi + \tau_2 \wedge \varphi.
\end{array}
\end{equation} 
The differential forms $\tau_0,\tau_1,\tau_2,\tau_3$ are called {\em intrinsic torsion forms} of the $\G_2$-structure $\f$, and they can be identified  with 
the components of the intrinsic torsion $T_\f$ belonging to the $\G_2$-modules $\mathcal{X}_1,\mathcal{X}_4,\mathcal{X}_2,\mathcal{X}_3$, respectively. 

Some classes of $\G_2$-structures with the defining conditions are recalled in Table \ref{G2classes}.

\begin{table}[ht]
\centering
\renewcommand\arraystretch{1.1}
\begin{tabular}{|c|c|c|}
\hline
class & type & conditions\\ \hline \hline
${\mathcal X}_1$          				& nearly parallel 			& $\tau_1, \tau_2, \tau_3 =0$  	\\ \hline
${\mathcal X}_2$               			& closed, calibrated	  		& $\tau_0, \tau_1, \tau_3 =0$	\\ \hline
${\mathcal X}_4$               			& locally conformal parallel 	& $\tau_0, \tau_2, \tau_3 =0$  	\\ \hline
${\mathcal X}_1\oplus {\mathcal X}_3$  	& coclosed, cocalibrated		& $\tau_1, \tau_2 = 0$ 		\\ \hline
${\mathcal X}_2 \oplus {\mathcal X}_4$     & locally conformal calibrated    & $\tau_0, \tau_3 =0$ 		\\ \hline
\end{tabular}
\vspace{0.1cm}
\caption{Some classes of $\G_2$-structures}\label{G2classes}
\end{table}
\renewcommand\arraystretch{1}

\subsection{Link with SU(3)-structures}\label{sectSU3}
An $\SU(3)$-structure on a six-dimensional manifold $N$ is the data of an almost Hermitian structure $(g,J)$ with fundamental 2-form $\omega\coloneqq g(J\cdot,\cdot)$ 
and a complex volume form $\Psi=\psi+\mathrm{i}\,\widehat\psi\in\Omega^{3,0}(M)$ of nonzero constant length.  

By \cite{Hit1}, an $\SU(3)$-structure $(g,J,\Psi)$ is completely determined by the real 2-form $\omega$ and the real 3-form $\psi$. 

Since $\G_2$ acts transitively on the 6-sphere with isotropy $\SU(3)$, every $\G_2$-structure on a 7-manifold $M$ induces an $\SU(3)$-structure on each oriented hypersurface. 
In particular, if $M$ is endowed with a torsion-free $\G_2$-structure $\f$, and $N\subset M$ is an oriented hypersurface, then $\f$ induces an $\SU(3)$-structure $(\omega,\psi)$ 
on $N$ which is {\em half-flat} according to the definition given in \cite{ChSa}. This means that the differential forms $\omega$ and $\psi$ satisfy the conditions 
\[
d(\omega\W\omega)=0,\quad d\psi=0.
\]
The inverse problem, i.e., establishing whether a half-flat $\SU(3)$-structure on a 6-manifold is induced by an immersion into a 7-manifold 
with a torsion-free $\G_2$-structure, can be analyzed using the so-called {\em Hitchin flow equations} (see \cite{Bry2,Hit1} for details). 

We now recall the definition of some special types of half-flat $\SU(3)$-structures. 
\begin{definition}\label{CpdSHF}
A half-flat $\SU(3)$-structure $(\omega,\psi)$ such that $d\omega=c\psi$ for some real number $c$ is said to be {\em coupled} if $c\neq0$, 
while it is called {\em symplectic half-flat} if $c=0$, i.e., if the 2-form $\omega$ is symplectic. 
A coupled $\SU(3)$-structure satisfying the additional condition $d\widehat\psi = -\frac{2}{3}c\,\omega\W\omega$ is called {\em nearly K\"ahler}. 
\end{definition}

If $N$ is a 6-manifold endowed with an $\SU(3)$-structure $(\omega,\psi)$, then the product manifold $N\times \R$ admits a $\G_2$-structure defined by the 3-form 
\[
\f\coloneqq \omega\W dt + \psi,
\] 
where $dt$ is the global 1-form on $\R$. Such $\f$ induces the product metric $g_\f=g+dt^2$. 
Moreover, $\varphi$  is calibrated (resp.~locally conformal calibrated)  if the  $\SU(3)$-structure $(\omega, \psi)$  is  symplectic half-flat   (resp.~coupled),   
while $\varphi$ is locally conformal parallel  if $(\omega, \psi)$  is nearly K\"ahler.

%%%%%%%%%%%%%%%%%%%%%%%%%%%%%%%%%%%%%%%%%%%%%%%%%%%%%%%%%%%%%%%%%%%%%%%%%%%%%%%%%%%%%%%%%
%%%%%%%%%%%%%%%%%%%%%%%%%%%%%%%%%%%%%%%%%%%%%%%%%%%%%%%%%%%%%%%%%%%%%%%%%%%%%%%%%%%%%%%%%
%																SPECIAL METRICS 
%%%%%%%%%%%%%%%%%%%%%%%%%%%%%%%%%%%%%%%%%%%%%%%%%%%%%%%%%%%%%%%%%%%%%%%%%%%%%%%%%%%%%%%%%
%%%%%%%%%%%%%%%%%%%%%%%%%%%%%%%%%%%%%%%%%%%%%%%%%%%%%%%%%%%%%%%%%%%%%%%%%%%%%%%%%%%%%%%%%
\section{$\G_2$-structures and special  metrics}\label{G2SpecMetrSect}
By \cite{Bry}, the Ricci tensor and the scalar curvature of the metric induced by a $\G_2$-structure $\f$ can be expressed
in terms of  the intrinsic torsion forms  $\tau_i$. In particular, the scalar curvature is given by
\[
\Scal(g_\f) = 12d^*\tau_1 + \frac{21}{8}\tau_0^2 + 30|\tau_1|^2-\frac{1}{2}|\tau_2|^2-\frac{1}{2}|\tau_3|^2,
\]
where $|\cdot|$ denotes the pointwise norm induced by $g_\f$.
Consequently, it has a definite sign for certain classes of $\G_2$-structures. 
For instance, when $\f$ is calibrated, then $\Scal (g_\f) = -\frac12 |\tau_2|^2$ is non-positive, while  
a nearly-parallel $\G_2$-structure always induces an Einstein metric with  positive scalar curvature $\Scal(g_\f)  =  \frac{21}{8} \tau_0^2.$

A  generalization of Einstein metrics is given by Ricci solitons. We recall the definition here. 
\begin{definition}\label{RicSol} 
A (complete) Riemannian metric $g$ on a smooth manifold $M$  is a 
{\em Ricci soliton} if its Ricci tensor satisfies the  equation
\[
\Ric(g) =  \lambda g + \mathcal{L}_X g,
\]
for some real constant $\lambda$ and some (complete) vector field $X$, where $\mathcal{L}$  denotes the Lie derivative. 
If in addition $X$ is the gradient of a smooth function $f \in\mathcal{C}^\infty(M)$, i.e., $X  = \nabla f$, 
then $g$  is said to be of {\em gradient type}.
\end{definition}

Equivalently, a Riemannian metric $g$ on $M$ is a Ricci soliton if  and only if there exists a  positive real valued function 
$h(t)$ and a family of diffeomorphisms $\eta_t$ such that $g(t) = h(t)\, \eta_t^* (g)$ is a solution of the Ricci flow starting from $g$ 
(see e.g.~\cite[Lemma 2.4]{ChKn}). 

Depending on the sign of $\lambda$, a Ricci solitons is called {\em expanding} ($\lambda<0$),  
{\em steady} ($\lambda=0$) or {\em shrinking} ($\lambda>0$).
Moreover, a Ricci soliton is said to be {\em trivial} if it is either Einstein or the product of a homogeneous Einstein metric with the Euclidean metric.
According to \cite{Ive}, if $M$ is a compact manifold with a Ricci soliton $g$ which is steady or expanding, then $g$ is Einstein.

A special class of Ricci solitons is given by {\em homogeneous} ones, which are defined as follows 
\begin{definition} 
A Ricci soliton $g$ on a smooth manifold $M$ is {\em homogeneous} if its isometry group acts transitively on $M.$ 
\end{definition}

 Properties of non-trivial homogeneous Ricci solitons were given by Lauret in \cite{Lau00}. In particular, he proved the following.
\begin{proposition}[\cite{Lau00}] 
Let $g$ be a non-trivial homogeneous Ricci soliton on a smooth manifold $M$. 
Then, $g$ is expanding and it cannot be of gradient type. Moreover, $M$ has to be non-compact.
\end{proposition} 

Currently, all known examples of nontrivial homogeneous Ricci solitons are {\em solvsolitons}, that is left-invariant Ricci solitons on simply connected solvable Lie groups.

Since requiring that the metric induced by a $\G_2$-structure is Einstein might impose some constraints on the intrinsic torsion, 
a natural problem is to investigate which types of $\G_2$-structures can induce an Einstein (or, more generally, a Ricci soliton) non-Ricci-flat metric, 
and to see whether there is any difference between the compact and noncompact cases.
For instance, if $M$ is a 7-manifold endowed with a locally conformal nearly parallel $\G_2$-structure $\f$ (torsion class $\mathcal{X}_1\oplus\mathcal{X}_4$)  
with $g_{\varphi}$ complete and Einstein, then $(M, \varphi)$ is either nearly parallel or conformally equivalent to the standard $7$-sphere (\cite{ClIv2}).

In what follows, we consider the cases of calibrated and locally conformal calibrated $\G_2$-structures. 

\subsection{Calibrated G$_{\mathbf2}$-structures}
A calibrated $\G_2$-structure $\f$ satisfies the equations 
\[
d\f =0,\quad d*_\f\f = \tau_2 \wedge \f,
\] 
with $\tau_2 \in \Omega_{14}^2(M)$.   
We collect some known properties of such type of $\G_2$-structures in the next results.

\begin{proposition}[\cite{Bry}] 
Let $\f$ be a calibrated $\G_2$-structure on $M$. Then,
\begin{enumerate}[1)]
\item $\Scal (g_\f) \leq 0$ and $\Scal (g_{\f}) =0$ if and only if $g_{\f}$ is Ricci-flat;
\item $\f$ defines an Einstein metric on $M$ if and only if  $d *_\f \f  = \tau_2 \wedge \f,$   
with $d \tau_2 =  \frac{3} {14} | \tau_2 |^2 \varphi + \frac 12 *_\f  (\tau_2 \wedge \tau_2).$
\end{enumerate}
\end{proposition}

\begin{corollary}[\cite{Bry, ClIv}]\label{CompactClosedEinstein} 
 Let $M$ be a compact 7-manifold with a calibrated $\G_2$-structure $\f$.
If the underlying metric $g_\f$ is Einstein, then $d *_\f  \f =0$ or, equivalenty, the holonomy group of 
$g_\f$ is a subgroup of $\G_2$.
\end{corollary}

The proof of the corollary follows from the identity  $d \left (  \frac 13 \tau_2^3  \right ) = \frac{2}{7} | \tau_2 |^4 *_\f1$ and Stokes' theorem. 
In detail, $\tau_2$ must vanish identically since
\[
0 = \int_{M} d \left (  \frac 13 \tau_2^3  \right ) = \int_{M}  \frac{2}{7} | \tau_2 |^4 *_\f  1.
\]

In the non-compact case, there is a known non-existence result involving the {\em $*$-Ricci tensor} and the 
{\em $*$-scalar curvature}, where
\[
\Ric^*(g_{\varphi})_{sm} = R_{ijkl}  \varphi_{ijs} \varphi_{klm}, \quad  \Scal^*(g_{\varphi}) = \mathrm{tr}_{g_{\varphi}} (\Ric^*(g_{\varphi})).
\]
 Such a result can be stated as follows.

\begin{theorem}[\cite{ClIv}]
Let $\f$ be a calibrated $\G_2$-structure on a 7-manifold $M.$  
If $g_{\f}$ is Einstein and  $*$-Einstein, i.e.,  $\Ric^* (g_{\f}) = \frac{\Scal^*(g_{\f}) }{7} g_{\f}$,  then  $g_{\f}$ is Ricci-flat.
\end{theorem}

In light of the previous results, one might investigate the existence of calibrated $\G_2$-structures that are Einstein but non-Ricci-flat on non-compact manifolds. 
This problem can be viewed as a $\G_2$-analogue of the Goldberg conjecture \cite{Gol}, which states that a compact Einstein 
almost-K\"ahler manifold has to be K\"ahler.   
Recall that a non-compact homogeneous example of Einstein strictly almost K\"ahler 6-manifold was constructed in  \cite{ADM}.

In the homogeneous setting, an answer to the above problem for calibrated $\G_2$-structures was given in \cite{FFM}.

All known examples of non-compact homogeneous Einstein manifolds are solvmanifolds, that is, simply connected solvable Lie groups 
endowed with a left-invariant Einstein metric.
The long-standing {\em Alekseevskii conjecture} \cite[Question 7.5]{Bes} states that a connected homogeneous Einstein space $\G/\K$ of negative scalar curvature 
must be diffeomorphic  to the Euclidean space.
Thus, Einstein solvmanifolds might exhaust the class of non-compact homogeneous Einstein manifolds. 
The conjecture is known to be true in dimensions five and lower by \cite{Jen,Nik}, and in dimension seven by \cite{ArLa}. 
So, seven-dimensional non-compact homogeneous Einstein manifolds are necessarily solvmanifolds. 

We now review some general results about Einstein metrics on solvmanifolds of arbitrary dimension. 
\begin{theorem}[\cite{Lau}] 
Every Einstein solvmanifold $(\mathrm{S}, g)$ is standard, i.e., the corresponding solvable metric Lie algebra $(\frs, \langle  \cdot, \cdot  \rangle)$ 
admits the orthogonal decomposition $\frs = \frn \oplus  \fra$, with $\frn =  [\frs, \frs]$ and $\fra$ abelian.
\end{theorem}
Recall that the dimension of the abelian summand $\fra$ in the decomposition 
\[
\frs = \frn \oplus  \fra
\]
is called the {\em rank} of the standard solvable metric Lie algebra $(\frs,\langle\cdot,\cdot\rangle)$. 

In contrast to the compact homogeneous case  (see e.g.~\cite[$\S$5]{Heb} and the references therein),  standard Einstein metrics are essentially unique.
\begin{theorem}[\cite{Heb}]  
A  standard Einstein metric is  unique up to isometry and scaling among invariant metrics.
\end{theorem}

\begin{remark}\  
\begin{enumerate} 
\item The study of standard Einstein solvmanifolds reduces to those with $\dim{\mathfrak a}  = 1$  
 (cf.~\cite[Thm.~4.18]{Heb}).
\item The Lie algebra of any standard Einstein solvmanifold  resembles an Iwasawa subalgebra of a semisimple
Lie algebra, since $\mathrm{ad}_A$ is symmetric and  non-zero for any $A \neq 0  \in \mathfrak a$, 
and there exists some $A^{\sst0}\in\fra$ such that $\mathrm{ad}_{A^{\sst0}}|_\frn$ is positive definite (see \cite[Thm.~4.10]{Heb}).
\end{enumerate}
\end{remark}

Using the rank  of a standard solvable metric Lie algebra, it is possible to get a classification of  seven-dimensional Einstein solvmanifolds 
(see e.g.~\cite[Thm.~4.4]{FFM}).  Then, using the
obstructions to the existence of calibrated $\G_2$-structures on Lie algebras given in  \cite{CoFe}, we have 
the following result.

\begin{theorem}[\cite{FFM}]\label{ClosedEinsteinSolvable}
Let  $g_{\varphi}$ be the metric determined by a left-invariant calibrated $\G_2$-structure $\f$ on a solvmanifold. 
Then,  $g_{\varphi}$ is Einstein if and only if $g_{\varphi}$ is flat.
\end{theorem}

\begin{remark} 
 Note that a similar theorem can be proved also for cocalibrated $\G_2$-structures \cite{FFM}.
Moreover, Theorem \ref{ClosedEinsteinSolvable} shows that left-invariant calibrated $\G_2$-structures 
behave differently from almost K\"ahler structures \cite{ADM}.
\end{remark}

The situation is different if we require that $g_{\varphi}$ is a non-trivial Ricci soliton.
Indeed, non-compact examples of manifolds admitting a calibrated $\G_2$-structure inducing 
a non-trivial  Ricci soliton were constructed in \cite{FFM}, and they are all  {\em nilsolitons} 
(see Theorem \ref{ClosedG2Nilsoliton} and Example \ref{FFMClosedRS} below).

\begin{definition} 
Let $\mathrm{N}$ be a simply connected nilpotent Lie group endowed with a left-invariant Riemannian metric $g$, and denote by 
$(\frn,\langle\cdot,\cdot\rangle)$ the corresponding metric nilpotent Lie algebra. 
The metric $g$ is called {\em nilsoliton} if its Ricci endomorphism $\Ric(g)$ on $\frn$ 
differs from a derivation $D$ of $\frn$ by a scalar multiple of the identity map $I$ , i.e., 
\[
\Ric(g) =  \lambda I + D,
\]
for some real number $\lambda$.  
\end{definition}

By \cite[Prop.~1.1]{Lau0}, a left-invariant Riemannian metric on a simply connected nilpotent Lie group is a nilsoliton if and only if 
it is a Ricci soliton according to Definition \ref{RicSol}. 
It is worth recalling here that non-abelian nilpotent Lie groups cannot admit any left-invariant Einstein metric unless it is flat \cite{Mil}.

\begin{remark}
As the existence of a nilsoliton on a simply connected nilpotent Lie group $\mathrm{N}$ implies the existence of a non-zero symmetric derivation on the corresponding 
nilpotent Lie algebra $\frn$, nilsolitons might not exist. This is the case, for instance, of Lie algebras having nilpotent derivation algebra. 
Such Lie algebras are nilpotent by Engel's Theorem, and they are known as {\em characteristically nilpotent} in literature.
\end{remark}

Before reviewing some properties of nilsolitons, we recall the following. 
\begin{definition}  
Let $({\mathfrak n}, [ \cdot , \cdot ]_{\mathfrak n}, \langle \cdot,\cdot \rangle_\frn)$ be a metric nilpotent Lie algebra. 
A metric Lie algebra $({\mathfrak s} = {\mathfrak n} \oplus {\mathfrak a}, [\cdot,\cdot],  \langle\cdot, \cdot \rangle)$ is a {\em metric solvable extension} 
of $({\frn},  [ \cdot , \cdot ]_{\mathfrak n},  \langle \cdot, \cdot \rangle_\frn)$ 
if  the restriction to $\frn$ of the Lie bracket $[\cdot , \cdot ]$  of $\mathfrak s$ coincides with $[ \cdot , \cdot ]_{\mathfrak n}$ 
and $\langle \cdot, \cdot \rangle \vert_{{\mathfrak n} \times {\mathfrak n}} = \langle \cdot, \cdot \rangle_\frn$.
\end{definition}
 
\begin{theorem}[\cite{Lau0}]\label{nilsolitonsLau} 
Let $\mathrm{N}$ be a simply connected nilpotent Lie group with Lie algebra $\frn$. Then,
\begin{enumerate}[1)]
\item A  nilsoliton metric  on $\mathrm{N}$  is unique up to isometry and scaling; 
\item $\mathrm{N}$ has a nilsoliton metric $g$ if and only if the corresponding metric Lie algebra $(\frn,\langle\cdot,\cdot\rangle)$ 
is an Einstein nilradical,  i.e., it has a metric solvable extension $\frs=\frn\oplus\fra$, with $\fra$ abelian, whose corresponding solvmanifold is Einstein. 
\end{enumerate}
\end{theorem}

 From now on, we will use the following notation to define a Lie algebra.
Suppose that $\frg$ is a seven-dimensional Lie algebra, whose dual space ${\frg}^*$ is spanned by $\{ e^1,\ldots ,e^7\}$ satisfying
$$
de^i=0, \quad 1\leq i\leq 4, \qquad de^5=e^{12}, \qquad de^6=e^{13}, \qquad de^7=0,
$$
where $d$ is the Chevalley-Eilenberg differential of $\frg$. Then, we will write
$$
\frg \,=\,(0, 0, 0, 0, e^{12}, e^{13}, 0)
$$
with the same meaning. 

In order to show the existence of nilpotent Lie algebras with a calibrated $\G_2$-structure inducing a nilsoliton,
we need to recall the classification of the nilpotent Lie algebras admitting a calibrated $\G_2$-structure given in \cite{CoFe}.

\begin{theorem}[\cite{CoFe}]\label{CoFeClass}
Up to isomorphism, there are exactly twelve nilpotent Lie algebras admitting a calibrated $\G_2$-structure. They are:
\[
\renewcommand\arraystretch{1.2}
\begin{array}{rcl}
\frn_1	&=&	(0, 0, 0, 0, 0, 0, 0),\\
\frn_2 	&=& (0, 0, 0, 0, e^{12}, e^{13}, 0),\\
\frn_3 	&=& (0, 0, 0, e^{12}, e^{13}, e^{23}, 0),\\
\frn_4 	&=& (0, 0, e^{12}, 0, 0, e^{13} + e^{24}, e^{15}),\\
\frn_5	&=& (0, 0, e^{12}, 0, 0, e^{13}, e^{14} + e^{25}),\\
\frn_6 	&=& (0, 0, 0, e^{12}, e^{13}, e^{14}, e^{15}),\\
\frn_7 	&=& (0, 0, 0, e^{12}, e^{13}, e^{14} + e^{23}, e^{15}),\\
\frn_8 	&=& (0, 0, e^{12}, e^{13}, e^{23}, e^{15} + e^{24}, e^{16} + e^{34}),\\
\frn_9 	&=& (0, 0, e^{12}, e^{13}, e^{23}, e^{15} + e^{24}, e^{16} + e^{34} + e^{25}),\\
\frn_{10}	&=& (0, 0, e^{12}, 0, e^{13} + e^{24}, e^{14}, e^{46} + e^{34} + e^{15} + e^{23}),\\
\frn_{11} 	&=& (0, 0, e^{12}, 0, e^{13}, e^{24} + e^{23}, e^{25} + e^{34} + e^{15} + e^{16} - 3 e^{26}),\\
\frn_{12} 	&=& (0, 0, 0, e^{12}, e^{23},-e^{13}, 2 e^{26} - 2 e^{34} - 2 e^{16} + 2 e^{25}).
\end{array}
\renewcommand\arraystretch{1}
\]
\end{theorem}

\smallskip

Comparing the previous classification with the results in \cite{Fer},   it turns out that, up to  isomorphism,   
$\frn_9$ is the unique nilpotent Lie algebra with a calibrated $\G_2$-structure but not admitting any nilsoliton. Moreover, 
the existence of a nilsoliton on the Lie algebra $\frn_{10}$ was shown in \cite[Example 2]{Fer1}, but its explicit expression is not known. 
Therefore, it is still an open problem to determine whether the Lie algebra $\frn_{10}$ admits a calibrated $\G_2$-structure inducing the nilsoliton.
For the remaining Lie algebras, we have the following. 
\begin{theorem}[\cite{FFM2}]\label{ClosedG2Nilsoliton}
Up to isomorphism,  $\frn_2, \frn_4, \frn_6$ and $\frn_{12}$  are the unique s-step
nilpotent Lie algebras ($s = 2, 3$) with a nilsoliton inner product determined by a calibrated $\G_2$-structure. 
\end{theorem}

\begin{remark}
Note that the Lie algebra $\frn_i$, $i=3, 5, 7, 8,11,$ has a nilsoliton inner
product but no calibrated $\G_2$-structure defining the nilsoliton \cite{FFM2}.
\end{remark}

In the next example, we write the expression of a calibrated $\G_2$-structure inducing the nilsoliton inner product on $\frn_i$, for $i=2,4,6,12$. 
Moreover, in each case we also specify the negative real number $\lambda$ and the derivation $D$ of $\frn_i$ for which 
 $\Ric=\lambda I + D$.

\begin{example}[\cite{FFM2}]\label{FFMClosedRS}
Consider the nilpotent Lie algebras $\frn_2, \frn_4, \frn_6$ with the structure equations given in 
Theorem \ref{CoFeClass}. Then,
\[
\renewcommand\arraystretch{1.5}
\begin{array}{rl}
\frn_2:		&  \f_{\sst2}   = e^{147}+e^{267}+e^{357}+e^{123}+e^{156}+e^{245}-e^{346},	\\
			&  \lambda=-2,\,\,\,\, ~D=\diag\left(1, \frac 32, \frac 32, 2, \frac{5}{2}, \frac{5}{2}, 2 \right); \\
\frn_4:		&  \f_{\sst4}   =	-e^{124}-e^{456}+e^{347}+e^{135}+e^{167}+e^{257}-e^{236},	\\
			&  \lambda=-\frac52,\,\,\,\, ~D=\diag\left(1, \frac 32, \frac 52, 2, 2, \frac 72, 3 \right);	\\
\frn_6:		&  \f_{\sst6}   = e^{123}+e^{145}+e^{167}+e^{257}-e^{246}+e^{347}+e^{356},	\\
			&  \lambda=-\frac52,\,\,\,\, ~D=\diag\left( \frac 12, 2, 2, \frac 52, \frac 52, 3, 3 \right ).	
\end{array}
\]
For $\frn_{12}$, we firstly consider a basis $\{e^1,\ldots,e^7\}$ of its dual  space
${\frn_{12}}^*$ for which the structure equations are
\[
\left(0,0,0,\frac{\sqrt{3}}{6}e^{12}, \frac{\sqrt{3}}{12}e^{13}-\frac{1}{4}e^{23}, -\frac{\sqrt{3}}{12}e^{23}-\frac{1}{4}e^{13},
\frac{\sqrt{3}}{12}e^{16}-\frac{\sqrt{3}}{6}e^{34}+\frac{\sqrt{3}}{12}e^{25}+\frac{1}{4}e^{26}-\frac{1}{4}e^{15}\right).
\]
Then, a calibrated $\G_2$-structure satisfying the required  properties is 
\[
\f_{\sst12} = -e^{124}+e^{167}+e^{257}+e^{347}-e^{456}+ e^{135}-e^{236}, 
\]
with $\lambda=-\frac14$  and $D=\frac18\,\diag(1,1,1,2,2,2,3)$. 
\end{example}

\begin{remark}
The nilsoliton condition is less restrictive for cocalibrated $\G_2$-structures. Indeed, on each 2-step nilpotent Lie algebra admitting cocalibrated $\G_2$-structures there 
exists a cocalibrated $\G_2$-structure inducing the nilsoliton inner product (see \cite{BFF}). 
\end{remark}

\subsection{Locally conformal calibrated G$_{\mathbf2}$-structures}
A $\G_2$-structure $\f$ is said to be {\em locally conformal calibrated} if the intrinsic torsion forms $\tau_0$ and $\tau_3$ vanish identically (cf.~Table \ref{G2classes}). 
In this case, equations \eqref{IntTorFor} reduce to 
\[
d\f=3\tau_1\W\f,\quad d*_\f\f = 4\tau_1\W*_\f\f +\tau_2\W\f. 
\]
Let $\theta\coloneqq3\tau_1=-\frac{1}{4}*_\f(*_\f d\f\W\f)$ denote the {\em Lee form} of the $\G_2$-structure. 
Taking the exterior derivative of both sides of the equation $d\f=\theta\W\f$, we get $d\theta\W\f=0$. This implies $d\theta=0$. 
Consequently, each point of the manifold has an open neighborhood $\mathcal{U}$ where $\theta=df$ for some $f\in\mathcal{C}^\infty(\mathcal{U})$, and the 
3-form $e^{-f}\f$ defines a calibrated $\G_2$-structure on $\mathcal{U}$. 
Hence, locally conformal calibrated $\G_2$-structures are locally conformal equivalent to calibrated $\G_2$-structures. 

Motivated by Corollary \ref{CompactClosedEinstein} and Theorem \ref{ClosedEinsteinSolvable},  
it is natural to investigate the existence of locally conformal calibrated $\G_2$-structures whose associated metric is Einstein and non-Ricci-flat.  
In what follows, we refer to a locally conformal calibrated $\G_2$-structure $\f$ with $g_\f$ Einstein as an {\em Einstein locally conformal calibrated $\G_2$-structure}. 

On compact manifolds, the following constraint on the scalar curvature holds. 
\begin{theorem}[\cite{FiRa1}] 
An Einstein locally conformal calibrated $\G_2$-structure on a compact seven-dimensional manifold has non-positive scalar curvature. 
\end{theorem}

Since homogeneous Einstein manifolds with negative scalar curvature are non-compact (cf.~\cite[Thm.~7.4]{Bes}) and since every homogeneous Ricci-flat metric is flat (see \cite{AlKi}), 
an immediate consequence of the previous result is the following. 
\begin{corollary}[\cite{FiRa1}]\label{corFiRa1}
A compact homogeneous 7-manifold cannot admit an invariant Einstein locally conformal calibrated $\G_2$-structure $\f$ unless the underlying metric $g_\f$ is flat.
\end{corollary}

\begin{remark}
By \cite[Thm.~3.1]{Boh}, the result of Corollary \ref{corFiRa1} is valid more generally on every compact locally homogeneous space. 
\end{remark}

In the non-compact setting, there is an example of a simply connected solvable Lie group endowed with a left-invariant locally conformal calibrated $\G_2$-structure $\f$ 
such that $g_\f$ is Einstein non-Ricci-flat. 
Thus, the result of \cite{FFM} recalled in Theorem \ref{ClosedEinsteinSolvable} is not true anymore for locally conformal calibrated $\G_2$-structures. 
Before describing the example, we recall some useful results. 

Consider a 6-manifold $N$ endowed with a coupled $\SU(3)$-structure $(\omega,\psi)$, with $d\omega=c\psi$ (cf.~Definition \ref{CpdSHF}). 
As we mentioned in $\S$\ref{sectSU3}, the 3-form  $\f  \coloneqq \omega \wedge dt + \psi$  defines a  locally conformal calibrated $\G_2$-structure 
on the product manifold $N \times \R$. It is not difficult to check that the corresponding Lee form is $\theta=-c\,dt$.  

In \cite{FiRa1}, the classification of six-dimensional nilpotent Lie algebras admitting a coupled $\SU(3)$-structure inducing a nilsoliton was achieved (see also \cite[$\S$4.1]{FiRaSUSY}).  
We recall it in the next theorem.  

\begin{theorem}[\cite{FiRa1}] 
A non-abelian six-dimensional nilpotent Lie algebra admitting a coupled  $\SU(3)$-structure  is isomorphic to one of the following 
\[
\frh_1 = (0,0,0,e^{12},e^{14}-e^{23},e^{15}+e^{34}),\qquad \frh_{2}=(0,0,0,0,e^{13}-e^{24},e^{14}+e^{23}).
\] 
Moreover, the only one  admitting  a  coupled  $\SU(3)$-structure inducing a nilsoliton is $\frh_2$.
\end{theorem}

\begin{remark}
Notice that $\frh_2$ is the Lie algebra of the three-dimensional complex Heisenberg group.
\end{remark}

We are now ready to describe the example. 
\begin{example}[\cite{FiRa1}]
Consider the coupled $\SU(3)$-structure on $\frh_2$ defined by the pair
\[
\omega=e^{12}+e^{34}-e^{56},\quad \psi = e^{136}-e^{145}-e^{235}-e^{246}. 
\]
It satisfies the equation $d\omega=-\psi$, and it induces the nilsoliton inner product $g=\sum_{k=1}^6(e^k)^2$ with Ricci operator 
\[
\Ric=-3I +4 \, {\mbox{diag}} \left(\frac12,\frac12,\frac12,\frac12,1,1\right),
\]
where $D= {\rm diag}\left(\frac12,\frac12,\frac12,\frac12,1,1\right)$ is a symmetric derivation of $\frh_2$. 
Consequently, the metric rank-one solvable extension $\frs = \frh_2 \oplus  \langle e_7\rangle$ of $\frh_2$ with structure equations
\[
\left(\frac12e^{17},\frac12e^{27},\frac12e^{37},\frac12e^{47}, e^{13} - e^{24}+e^{57},  e^{14} + e^{23}+e^{67},0\right)
\]
is endowed with the Einstein (non-Ricci-flat) inner product $g + (e^7)^2.$ This is precisely the inner product $g_\f$ induced by the 3-form 
\[
\f=\omega\W e^7+\psi,
\]
which defines a locally conformal calibrated $\G_2$-structure on $\frs$. 
A simple computation shows that the non-vanishing intrinsic torsion forms of $\f$ are 
\[
\tu = - \frac{1}{3}\,e^7,\quad  \td=  -\left( \frac{5}{3}\,e^{12} + \frac{5}{3}\,e^{34} + \frac{10}{3}\,e^{56} \right).
\]

Clearly, left multiplication allows to extend $\f$ to a left-invariant Einstein locally conformal calibrated $\G_2$-structure 
on the simply connected nilpotent Lie group corresponding to $\frs$. 
\end{example}

We conclude this section recalling a general structure result for compact 7-manifolds endowed with a locally conformal calibrated $\G_2$-structure with nowhere 
vanishing Lee form. 
\begin{theorem}[\cite{FFR}]\label{StructResLCC}
Let $M$ be a compact, connected seven-dimensional manifold endowed with a locally conformal calibrated $\G_2$-structure $\f$, with nowhere vanishing Lee 
form $\theta.$
Suppose  that  $\mathcal{L}_X\f =0$, where $X$ is the $g_\f$-dual vector field of $\theta$. Then, 
\begin{enumerate}[1)]
\item  $M$ is the total space of a fibre bundle over $\mathbb{S}^1$, and each fibre is endowed with a coupled $\SU(3)$-structure; 
\item  $M$ has a locally conformal calibrated $\G_2$-structure $\hat \varphi$ such that $d \hat \varphi = \hat \theta \wedge \hat\varphi$,  
where $\hat \theta$ is a 1-form with integral periods.
\end{enumerate}
\end{theorem}

The previous theorem implies in particular that $M$ is the mapping torus of a diffeomorphism $\nu$ of a certain 6-manifold $N,$ i.e., 
$M$ is diffeomorphic to the quotient of 
$N\times \R$ by the infinite cyclic group of diffeomorphisms generated by $(p,t)\mapsto(\nu(p),t+1).$

\begin{remark}  
It is worth recalling here that compact locally conformal parallel $\G_2$-manifolds can be characterized as fibre bundles over $\mathbb{S}^1$ with compact nearly K\"ahler fibre 
(see \cite{IPP, Ver}).
\end{remark}

%%%%%%%%%%%%%%%%%%%%%%%%%%%%%%%%%%%%%%%%%%%%%%%%%%%%%%%%%%%%%%%%%%%%%%%%%%%%%%%%%%%%%%%%%
%%%%%%%%%%%%%%%%%%%%%%%%%%%%%%%%%%%%%%%%%%%%%%%%%%%%%%%%%%%%%%%%%%%%%%%%%%%%%%%%%%%%%%%%%
%																LAPLACIAN FLOW 
%%%%%%%%%%%%%%%%%%%%%%%%%%%%%%%%%%%%%%%%%%%%%%%%%%%%%%%%%%%%%%%%%%%%%%%%%%%%%%%%%%%%%%%%%
%%%%%%%%%%%%%%%%%%%%%%%%%%%%%%%%%%%%%%%%%%%%%%%%%%%%%%%%%%%%%%%%%%%%%%%%%%%%%%%%%%%%%%%%%
\section{The Laplacian flow on Lie groups}\label{LapFlowSect}
Consider a 7-manifold $M$ endowed with a calibrated $\G_2$-structure $\fz$. The {\em Laplacian flow} starting from $\fz$ is the initial value problem 
\begin{equation}\label{LapFlow}
\begin{cases}
\ddt \f(t) = \Delta_{\f(t)}\f(t),\\
d \f(t)=0,\\
\f(0)=\fz,
\end{cases}
\end{equation}
where $\Delta_{\f}$ denotes the Hodge Laplacian of the Riemannian metric $g_{\f}$ induced by $\f$. 
This flow was introduced by Bryant in \cite{Bry} to study seven-dimensional manifolds admitting calibrated $\G_2$-structures. 
Notice that the stationary points of the flow equation in \eqref{LapFlow} are harmonic $\G_2$-structures, which coincide with torsion-free $\G_2$-structures on compact manifolds. 

Short-time existence and uniqueness of the solution of \eqref{LapFlow} when $M$ is compact were proved in \cite{BrXu}. 
\begin{theorem}[\cite{BrXu}] 
Assume that  $M$ is compact. Then, the Laplacian flow \eqref{LapFlow} has a unique solution defined for a short time
$t\in[0,\varepsilon)$, with $\varepsilon$  depending on $\fz$. 
\end{theorem}

As a consequence of the condition $d\f(t)=0$, the solution $\f(t)$ must belong to the open set  
\[
[\fz]_{\sst+}\coloneqq[\fz]\cap\Omega^3_{\sst+}(M)
\]
in the cohomology class $[\fz]$ as long as it exists.  

\begin{remark}
By \cite{Bry,Hit0}, the evolution equation in \eqref{LapFlow} is the gradient flow of Hitchin's volume functional 
\[
[\fz]_{\sst+}\ni \f \mapsto  \int_{M}dV_\f,
\]
with respect to a suitable $L^2$-metric on $[\fz]_{\sst+}$. 
\end{remark}

\subsection{Solutions to the Laplacian flow on nilpotent Lie groups}
Lie groups admitting left-invariant calibrated $\G_2$-structures constitute a convenient setting where it is possible to investigate the behaviour of the Laplacian flow in the 
non-compact case. 
In literature, results in this direction have been obtained on nilpotent and solvable Lie groups in various works \cite{FFM2, FiRa, Lau2, Lau3,Nic}. 
In the non-solvable case, the first examples of calibrated $\G_2$-structures have been exhibited only recently in \cite{FiRa3}, 
and the study of the Laplacian flow starting from any of them is still in progress. 

The main peculiarity of the known non-compact examples is that the solution of \eqref{LapFlow} exists on an infinite time interval.  
We recall here a result of \cite{FFM2}, while we refer the reader to \cite{Lau2,Lau3} for further examples. 

\begin{example}[\cite{FFM2}]\ 
\begin{enumerate}[1)]
\item On the nilpotent Lie algebra $\frn_2$, the solution of the Laplacian flow starting from the calibrated $\G_2$-structure $\f_{\sst2}$ given in Example \ref{FFMClosedRS} is
\[
\f(t) = e^{147}+e^{267}+e^{357}+\Big(\frac{10}{3} t +1\Big)^{3/5}e^{123}+e^{156}+e^{245}-e^{346},
\]
where $t\in \left (-\frac{3}{10},+ \infty \right).$
\item On the nilpotent Lie algebra $\frn_{12}$, the solution of the Laplacian flow starting from the calibrated $\G_2$-structure $\f_{\sst12}$ given in Example \ref{FFMClosedRS} is
 \[
\f(t)= -e^{124}+e^{167}+e^{257}+e^{347}-e^{456}+ \left(\frac13t+1\right)^{3/4}(e^{135}-e^{236}),
\]
with $t \in (-3, +\infty)$.
\end{enumerate}
\end{example}

In the previous example, both the calibrated $\G_2$-structures considered as initial value for the Laplacian flow induce the nilsoliton inner product on 
the corresponding nilpotent Lie algebra (cf.~Theorem \ref{ClosedG2Nilsoliton} and Example \ref{FFMClosedRS}). 
Furthermore, using suitable analytic techniques it is possible to show that the solution $\f(t)$ of \eqref{LapFlow} with $\f(0)=\f_{\sst4}$ on $\frn_4$ 
and $\f(0)=\f_{\sst6}$ on $\frn_6$ exists for $t\in(T,+\infty)$ with $T<0$, see \cite[Thm.~4.7, Thm.~4.8]{FFM2}. 
In all cases, it is then possible to analyze the behaviour of the solution $\f(t)$ when $t\rightarrow+\infty$. 

\begin{theorem}[\cite{FFM2}]  
On the simply-connected nilpotent Lie groups $\mathrm{N}_i$, $i = 2,4,6,12$, the Laplacian flow starting from the left-invariant calibrated $\G_2$-structure $\f_{i}$ 
has a global solution defined for $t \in (T, + \infty)$, with $T<0$. 
Moreover, all solutions converge to a flat $\G_2$-structure when $t\rightarrow+\infty$. 
\end{theorem}

\begin{remark}
The nilpotent Lie algebra $\frn_2$ may be seen as a product algebra $\frn_2 = \frn'\oplus\R$ with $\dim(\frn')=6$ and $\R=\langle e_7\rangle$, 
and the calibrated $\G_2$-structure $\f_{\sst2}$ on it can be written as 
\[
\f_{\sst2}= \omega\W e^7+\psi, 
\]
where $(\omega,\psi)$ is a symplectic half-flat $\SU(3)$-structure on $\frn'$ (cf.~Definition \ref{CpdSHF}). 
Moreover, the solution of the Laplacian flow starting from $\f_{\sst2}$ at $t=0$ is of the form
\[
\f(t) = f(t)\,\omega(t)\W e^7+\psi(t),
\]
where $(\omega(t),\psi(t))$ is a family of symplectic half-flat $\SU(3)$-structures on $\frn'$, and the function 
$f:\left(-\frac{3}{10},+\infty\right)\rightarrow\R^{\sst+}$
is given by $f(t)=\left(\frac{10}{3}t+1\right)^{-1/10}$.  

This fact is a consequence of a more general result which holds for a suitable class of symplectic half-flat $\SU(3)$-structures and allows to construct new examples 
of solutions of \eqref{LapFlow} on solvable Lie groups. For more details, we refer the reader to \cite{FiRa}. 
\end{remark}

\begin{remark}
The investigation reviewed in this section can be carried out also for the {\em Laplacian coflow} for cocalibrated $\G_2$-structures \cite{KMT} 
and its modified version introduced in \cite{Gri}. It turns out that the behaviour of these flows on solvable Lie groups is slightly different  from the behaviour of the Laplacian flow. 
We refer the reader to \cite{BFF2,BaFi} for a detailed treatment. 
\end{remark}
\bigskip 

\noindent
{\bf Acknowledgements.} 
The first author was partially supported by MINECO-FEDER Grant MTM2014-54804-P
and Gobierno Vasco Grant IT1094-16, Spain. 
The second and third authors
were partially supported by GNSAGA of INdAM. 
The second author was also supported by PRIN 2015 
``Real and complex manifolds: geometry, topology and harmonic analysis'' of MIUR. 
All authors would like to thank the organizers of the workshop ``$\G_2$-manifolds and related topics'' (Fields Institute, Toronto, August 2017).

%%%%%%%%%%%%%%%%%%%%%%%%%%%%%%%%%%%%%%%%%%%%%%%%%%%%%%%%%%%%%%%%%%%%%%%%%%%%%%%%%%%%%%%%%
%%%%%%%%%%%%%%%%%%%%%%%%%%%%%%%%%%%%%%%%%%%%%%%%%%%%%%%%%%%%%%%%%%%%%%%%%%%%%%%%%%%%%%%%%
%																BIBLIOGRAPHY
%%%%%%%%%%%%%%%%%%%%%%%%%%%%%%%%%%%%%%%%%%%%%%%%%%%%%%%%%%%%%%%%%%%%%%%%%%%%%%%%%%%%%%%%%
%%%%%%%%%%%%%%%%%%%%%%%%%%%%%%%%%%%%%%%%%%%%%%%%%%%%%%%%%%%%%%%%%%%%%%%%%%%%%%%%%%%%%%%%%


\begin{thebibliography}{10}

\bibitem{AlKi}  
D.~V.~Alekseevski{\u\i} and B.~N.~Kimel'fel'd.
\newblock Structure of homogeneous {R}iemannian spaces with zero {R}icci curvature.
\newblock {\em Funkcional. Anal. i Prilo\v zen.}, {\bf 9} (2), 5--11, 1975.

\bibitem{ADM}
V.~Apostolov, T.~Dr{\u{a}}ghici, and A.~Moroianu.
\newblock A splitting theorem for {K}{\"a}hler manifolds whose {R}icci tensors have constant eigenvalues.
\newblock {\em Internat. J. Math.}, {\bf 12} (7), 769--789, 2001.

\bibitem{ArLa}
R.~M. Arroyo and R.~A. Lafuente.
\newblock The {A}lekseevskii conjecture in low dimensions.
\newblock {\em Math. Ann.}, {\bf 367} (1-2), 283--309, 2017.

\bibitem{BFF}  
L.~Bagaglini, M.~Fern\'andez, and A.~ Fino.
\newblock Coclosed $\G_2$-structures inducing nilsolitons. 
\newblock {\em Forum Math.},  {\bf 30} (1), 109--128, 2018. 

\bibitem{BFF2}
L.~Bagaglini, M.~Fern\'andez, and A.~ Fino. 
\newblock Laplacian coflow on the 7-dimensional Heisenberg group.  
arXiv:1704.00295 [math.DG].

\bibitem{BaFi}
L.~Bagaglini and  A.~ Fino.
\newblock The Laplacian coflow on almost-abelian Lie groups.
\newblock  {\em Ann.~Mat.~Pura Appl.}, {\bf 197} (6), 1855--1873, 2018.

\bibitem{Bes}
A.~L.~Besse.
\newblock {\em Einstein manifolds}, vol.~10 of {\em Ergebnisse der Mathematik und ihrer Grenzgebiete}.
\newblock Springer-Verlag, Berlin, 1987.

\bibitem{Boh}
C.~B{\"o}hm.
\newblock On the long time behavior of homogeneous {R}icci flows.
\newblock {\em Comment. Math. Helv.}, {\bf 90} (3), 543--571, 2015.

\bibitem{Bon}
E.~Bonan.
\newblock Sur des vari{\'e}t{\'e}s riemanniennes {\`a} groupe d'holonomie {G}$_{2}$\ ou {S}pin{$(7)$}.
\newblock {\em C.~R.~Acad.~Sci.~Paris S{\'e}r.~A-B}, {\bf 262}, A127--A129, 1966.

\bibitem{Bry87}
R.~L. Bryant.
\newblock Metrics with exceptional holonomy.
\newblock {\em Ann. of Math.}, {\bf 126} (3), 525--576, 1987.

\bibitem{Bry}
R.~L.~Bryant.
\newblock Some remarks on {G$_2$}-structures.
\newblock In {\em Proceedings of {G}{\"o}kova {G}eometry-{T}opology {C}onference 2005},  pp. 75--109. G{\"o}kova Geometry/Topology Conference (GGT), G{\"o}kova, 2006.

\bibitem{Bry2}
R.~L. Bryant.
\newblock Non-embedding and non-extension results in special holonomy.
\newblock In {\em The many facets of geometry},  pp.  346--367. Oxford Univ. Press, Oxford, 2010.

\bibitem{BrXu}
R.~L.~Bryant and F.~Xu.
\newblock Laplacian flow for closed {G}$_2$-structures: Short time behavior.
arXiv:1101.2004 [math.DG]. 

\bibitem{ChSa}
S.~Chiossi and S.~Salamon.
\newblock The intrinsic torsion of {$\rm SU(3)$} and {G$_2$} structures.
\newblock In {\em Differential geometry, {V}alencia, 2001},  pp.  115--133. World Sci. Publ., River Edge, NJ, 2002.  

\bibitem{ChKn}
B.~Chow and D.~Knopf.
\newblock {\em The {R}icci flow: an introduction}, vol. 110 of {\em Mathematical Surveys and Monographs}.
\newblock Amer. Math. Soc., Providence, RI, 2004.

\bibitem{ClIv}
R.~Cleyton and S.~Ivanov.
\newblock On the geometry of closed {G}$_2$-structures.
\newblock {\em Comm. Math. Phys.}, {\bf 270} (1), 53--67, 2007.

\bibitem{ClIv2}  
R.~Cleyton and S.~Ivanov. 
\newblock Conformal equivalence between certain geometries in dimension 6 and 7. 
\newblock {\em Math. Res. Lett.}, {\bf 15} (4), 631--640, 2008.

\bibitem{CoFe}
D.~Conti and M.~Fern{{\'a}}ndez.
\newblock Nilmanifolds with a calibrated {G}$_2$-structure.
\newblock {\em Differ. Geom. Appl.}, {\bf 29} (4), 493--506, 2011.


\bibitem{Fer1} 
E. A.~Fern\'andez-Culma. Classification of 7-dimensional Einstein Nilradicals,
\newblock \emph{Transform. Groups}, {\bf 17} (3), 639--656, 2012.


\bibitem{Fer}  
E.~A.~Fern\'andez-Culma. 
\newblock Classification of Nilsoliton metrics in dimension seven. 
\newblock {\em J. Geom. Phys.}, {\bf 86}, 164--179, 2014.

\bibitem{FFM}
M.~Fern{{\'a}}ndez, A.~Fino, and V.~Manero.
\newblock {G}$_2$-structures on {E}instein solvmanifolds.
\newblock {\em Asian J. Math.}, {\bf 19} (2), 321--342, 2015.

\bibitem{FFM2} 
M.~Fern\'andez, A.~Fino, and V.~Manero. 
\newblock Laplacian flow of closed {G$_2$}-structures inducing nilsolitons.
\newblock {\em J. Geom. Anal.}, {\bf 26} (3), 1808--1837, 2016.   

\bibitem{FFR} 
M.~Fern\'andez, A.~Fino, and A.~Raffero.
\newblock Locally conformal calibrated {G}$_2$-manifolds
\newblock  {\em Ann.~Mat.~Pura Appl.}, {\bf 195} (5), 1721--1736, 2016.

\bibitem{FeGr}
M.~Fern{{\'a}}ndez and A.~Gray.
\newblock Riemannian manifolds with structure group {G}$_{2}$.
\newblock {\em Ann. Mat. Pura Appl.}, {\bf 132}, 19--45, 1982.

\bibitem{FiRaSUSY}
A.~Fino and A.~Raf{}fero. 
\newblock Coupled SU$(3)$-structures and supersymmetry. 
\newblock {\em Symmetry}, {\bf7} (2), 625--650, 2015.


\bibitem{FiRa1}
A.~Fino and A.~Raf{}fero.
\newblock Einstein locally conformal calibrated {${\rm G}_2$}-structures.
\newblock {\em Math. Z.}, {\bf 280} (3-4), 1093--1106, 2015.

\bibitem{FiRa} 
A.~Fino and A.~Raf{}fero. 
\newblock Closed warped $\G_2$-structures evolving under the Laplacian flow.  
\newblock To appear in {\em Ann.~Scuola Norm.~Sup.~Pisa Cl.~Sci.} doi: 10.2422/2036-2145.201709$\textunderscore$004.


\bibitem{FiRa3}
A.~Fino and A.~Raffero.
\newblock Closed G$_2$-structures on non-solvable Lie groups.
\newblock arXiv:1712.09664 [math.DG].

\bibitem{Gol}
S.~I.~Goldberg.
\newblock Integrability of almost {K}\"ahler manifolds.
\newblock {\em Proc. Amer. Math. Soc.}, {\bf 21}, 96--100, 1969.

\bibitem{Gri} 
S.~Grigorian. 
\newblock Short-time behavior of a modified Laplacian coflow of $\G_2$-structures. 
\newblock {\em Adv. Math.},  {\bf 248}, 378--415, 2013.

\bibitem{HaLa}
R.~Harvey and H.~B.~Lawson, Jr.
\newblock Calibrated geometries.
\newblock {\em Acta Math.}, {\bf 148}, 47--157, 1982.

\bibitem{Heb}
J.~Heber.
\newblock Noncompact homogeneous {E}instein spaces.
\newblock {\em Invent. Math.}, {\bf 133} (2), 279--352, 1998.


\bibitem{Hit0}
N.~Hitchin. 
\newblock The geometry of three-forms in six and seven dimensions.
\newblock {\em J. Differ. Geom.}, {\bf 55} (3), 547--576, 2000.

\bibitem{Hit1}
N.~Hitchin.
\newblock Stable forms and special metrics.
\newblock In {\em Global differential geometry: the mathematical legacy of {A}lfred {G}ray ({B}ilbao, 2000)}, vol.~288 of {\em Contemp. Math.}, 
 pp. 70--89. Amer. Math. Soc., 2001.

\bibitem{HWY}
H.~Huang, Y.~Wang and C.~Yao.  
\newblock Cohomogeneity-one $\G_2$-Laplacian flow on 7-torus. 
\newblock {\em J. London Math. Soc.}, {\bf 98} (2), 349--368, 2018. 


\bibitem{IPP}
S.~Ivanov, M.~Parton, and P.~Piccinni.
\newblock Locally conformal parallel {$\rm G_2$} and {${\rm Spin}(7)$} manifolds.
\newblock {\em Math. Res. Lett.}, {\bf 13} (2-3), 167--177, 2006. 

\bibitem{Ive}
T.~Ivey.
\newblock Ricci solitons on compact three-manifolds.
\newblock {\em Differ. Geom. Appl.}, {\bf 3} (4), 301--307, 1993.

\bibitem{Jen}
G.~R. Jensen.
\newblock Homogeneous {E}instein spaces of dimension four.
\newblock {\em J. Differ. Geom.}, {\bf 3}, 309--349, 1969.
  
\bibitem{KMT}
S.~Karigiannis, B.~McKay, and M.-P. Tsui. 
\newblock Soliton solutions for the {L}aplacian co-flow of some {$\G_2$}-structures with symmetry.
\newblock {\em Differ. Geom. Appl.}, {\bf 30} (4), 318--333, 2012.

\bibitem{Lau0}
J.~Lauret.
\newblock Ricci soliton homogeneous nilmanifolds.
\newblock {\em Math. Ann.}, {\bf 319} (4), 715--733, 2001.

\bibitem{Lau}
J.~Lauret.
\newblock Einstein solvmanifolds are standard.
\newblock {\em Ann. of Math.}, {\bf 172} (3), 1859--1877, 2010.

\bibitem{Lau00}
J.~Lauret.
\newblock Ricci soliton solvmanifolds
\newblock  {\em J.~Reine Angew.~Math.}, {\bf 650}, 1--21, 2011. 

\bibitem{Lau2}
J.~Lauret.
\newblock Laplacian flow of homogeneous {G$_2$}-structures and its solitons.
\newblock {\em Proc. Lond. Math. Soc.}, {\bf 114} (3), 527--560, 2017.

\bibitem{Lau3}
J.~Lauret.
\newblock Laplacian solitons: questions and homogeneous examples.
\newblock {\em Differ. Geom. Appl.}, {\bf 54} (B), 345--360, 2017.
 
\bibitem{LoWe1}
J.~D. Lotay and Y.~Wei.
\newblock Laplacian flow for closed {G$_2$} structures: {S}hi-type estimates, uniqueness and compactness.
\newblock {\em Geom. Funct. Anal.}, {\bf 27} (1), 165--233, 2017.

\bibitem{LoWe2}
J.~D. Lotay and Y.~Wei.
\newblock Stability of torsion-free G$_2$ structures along the Laplacian flow.
\newblock To appear in {\em J.~Differ. Geom.}, arXiv:1504.07771 [math.DG].

\bibitem{LoWe3}
J.~D. Lotay and Y.~Wei.
\newblock Laplacian flow for closed G$_2$ structures: real analyticity.
\newblock   To appear in {\em Comm.~Anal.~Geom.}, arXiv:1601.04258 [math.DG].
  
\bibitem{Mil}
J.~Milnor.
\newblock Curvatures of left invariant metrics on {L}ie groups.
\newblock {\em Adv. Math.}, {\bf 21} (3), 293--329, 1976.

\bibitem{Nic}
M.~Nicolini.
\newblock Laplacian solitons on nilpotent Lie groups.
\newblock To appear in  {\em B.~Belg.~Math.~Soc.}, arXiv:1608.08599 [math.DG].
 
\bibitem{Nik}
Y.~G. Nikonorov.
\newblock On the {R}icci curvature of homogeneous metrics on noncompact homogeneous spaces.
\newblock {\em Sibirsk. Mat. Zh.}, {\bf 41} (2), 421--429, iv, 2000.

\bibitem{Ver}
M.~Verbitsky.
\newblock An intrinsic volume functional on almost complex 6-manifolds and nearly {K}{\"a}hler geometry.
\newblock {\em Pacific J. Math.}, {\bf 235} (2), 323--344, 2008.

\end{thebibliography}
\end{document}